\documentclass[11pt]{article}

\usepackage{latexsym, amsmath}
\newtheorem{theorem}{Theorem}
\newtheorem{lemma}{Lemma}
\newtheorem{corollary}{Corollary}

\newenvironment{proof}{{\it Proof:\/}}{\hfill $\Box$\\ }

\newcommand{\F}{\mbox{\bf  F}}

\title{On the List and Bounded Distance Decodibility of the
Reed-Solomon Codes \\
(Extended Abstract)}

\date{}

\author{Qi Cheng\thanks{School of Computer Science,
the University of Oklahoma, Norman, OK 73019, USA. Email: {\tt
qcheng@cs.ou.edu}. This research is partially supported by NSF
Career Award CCR-0237845. } \and Daqing Wan\thanks{Department of
Mathematics, University of California, Irvine, CA 92697. Email:
{\tt dwan@math.uci.edu}. Institute of Mathematics, Chinese Academy
of Sciences, Beijing, P.R. China. Partially supported by NSF and
NSFC. } }

\begin{document}

\maketitle

\begin{abstract}
For an error-correcting code and a distance bound,
the list decoding problem is to compute all the codewords
within the given distance to a received message.
The {\em bounded distance decoding} problem, on the other hand,
is to find  one codeword if there
exists one or more codewords within the given distance,
or to  output the empty set if there does not.
Obviously the bounded distance decoding problem is not as hard
as the list decoding problem.
For a Reed-Solomon code $[n,k]_q$, a simple counting argument shows
that for any integer $g < n $, there exists at least one Hamming ball
of radius $n-g$, which contains at least $ { n \choose g} \over
q^{g-k} $ many codewords. Let $\hat{g}(n,k,q) $ be the smallest
integer $g$ such that $ {{ n \choose g} \over q^{g-k}} < 1 $.
For the distance bound between $n - \sqrt{nk}$ and
$n- \hat{g}(n,k,q) $, we do not know whether the Reed-Solomon code
is list, or bounded distance decodable, nor do we know whether there
are polynomially many codewords in all balls of the radius. It is
generally believed that the answers to both questions are no.  There
are public key cryptosystems proposed recently, whose security is
based on the assumptions.  In this paper, we prove: (1)
List decoding can not be done for radius $n -
\hat{g}(n,k,q) $ or larger, otherwise the discrete logarithm over
$\F_{q^{\hat{g}(n,k,q) -k}}$ is easy.  (2) Let $h$ be a positive
integer satisfying $ h < q^{1/4} -2 $.  We show that the
discrete logarithm problem over $\F_{q^{h}}$ can be efficiently
reduced to the bounded distance decoding problem of the Reed-Solomon
code $[q, 3h+4]_q$ with radius $q - 4h -4$.  These results
show that the decoding problems for the Reed-Solomon code are at
least as hard as the discrete logarithm problem over finite fields.
The main tools to obtain these results are an interesting connection
between the problems of list-decoding of
Reed-Solomon code and the problems of
discrete logarithms over finite fields, and a generalization
of the Katz's theorem,
which  concerns representations of elements in an extension
finite field by products of linear factors.
\end{abstract}

\section{Introduction and Motivation}


An error-correcting code $C$ over an alphabet $\Sigma$
is an injective map $\phi: \Sigma^k \rightarrow \Sigma^n$.
When we need to transmit a message of $k$ letters over a noisy channel,
we apply the map on the message
first ( i.e. encode the message ) and send its image (i.e. the codeword)
of $n$ letters  over the channel.
The Hamming distance between two sequence
of letters of the same length
is the number of positions where two sequences differ.
A good error-correcting  code should have a large {\em minimum
distance} $d$, which is defined to be the minimum Hamming distance
between any two codewords in $\phi(\Sigma^k)$.
A received message, possibly corrupted,
but with no more than $(d-1)/2$ errors, corresponds to
a unique codeword, thus may be decoded into the original
message despite errors occur during the communication.

Error-correcting codes are widely used in practice and are
mathematically interesting and intriguing.
It attracts the attention of theoretical
computer science community recently.  Several major achievements of
theoretical computer science, notably the Probabilistically
Checkable Proofs and derandomization techniques,
rely heavily  on the techniques in error-correcting codes.
We refer to  the survey \cite{Sudan01} for details.

For the purpose of efficient encoding and decoding,
$\Sigma$ is usually  set to be a finite field, and
the map $\phi$ is set to be linear.
Numerous error correcting codes  have been proposed,
among them, the
Reed-Solomon codes are particularly important.
They were deployed to transmit information from and
to spaceships, and were used to store information in optical media.
The Reed-Solomon code $[n,k]_q$,
is the map from $a_0, a_1, \cdots, a_{k-1} \in \F_q$
to $(a_0 + a_1x + \cdots + a_{k-1} x^{k-1})_{x \in S \subseteq \F_q}$
for some $|S| = n$.
(The choice of $S$ will not affect our results in this paper. )
Since any two different
polynomials with degree $k-1$ can share at most $k-1$
points, the minimum distance of the Reed-Solomon code
is $n-k+1$.
If the radius of a
Hamming ball is less than half of the minimum distance, there
should be at most one codeword in the Hamming ball. Finding
the codeword is called {\em unambiguous decoding}.
It was solved, see \cite{BerlekampWe86} for a simple algorithm.

If we gradually increase the
radius, there will be two or more codewords lying in some Hamming
balls. Can we  efficiently enumerate all the codewords in any
Hamming ball of certain radius? This is the so called list decoding
problem.  The notion was first introduced by Elias \cite{Elias57}.
There was virtually no  progress on this problem for radius
slightly larger than  half of the minimum distance,
until Sudan published his influential paper \cite{Sudan97}.
His result was subsequently improved,
the best algorithm \cite{GuruswamiSu99} solves the
list decoding problem for  radius as large as
$n - \sqrt{nk}$.
The work sheds new light on the limitation of list decodibility of Reed-Solomon
codes.  To the other extreme, if the radius
is greater than or equal to the minimum distance, there are
exponentially many codewords in some Hamming balls.


The decoding problem of Reed-Solomon codes can be formulated into the
problem of {\em curve fitting} or {\em polynomial reconstruction}.  In
the problem, we are given $n$ points $(x_1, y_1), (x_2, y_2), \cdots,
(x_n, y_n)$.  The goal is to find  polynomials of degree $k-1$ that
pass at least $g$ points.  In this paper, we only consider the case
when points have distinct $x$-coordinates.  If we allow multiple
occurrences of $x$-coordinates, the problem is NP-hard
\cite{GoldreichRu00}, and it is not relevant to the Reed-Solomon
decoding problem.  If $g\geq (n +k)/2$, it corresponds to the
unambiguous decoding of Reed-Solomon codes.  If $g > \sqrt{nk}$, the
radius is less than $n - \sqrt{nk}$, the problem can be solved by the
Guruswami-Sudan algorithm.  If $g \leq k$, it is possible that there
are exponentially many solutions, but finding one is very easy.

In this paper, we study the following question: How large
can we increase the radius before the list decoding problem or the
bounded distance decoding problem become infeasible?  The question
has been under intensive investigations for Reed-Solomon codes and other
error-correcting codes.  The case of general non-linear codes has been
solved \cite{GoldreichRu00}.  The case for linear codes is
much harder.  Some partial results have been obtained in
\cite{GuruswamiHa02, Guruswami02}. However,
none of them applies to Reed-Solomon codes.  No negative result is
known about the list decodibility of Reed-Solomon codes, except a
simple bound given by Justesen and Hoholdt \cite{JustesenHo01}, which
states that for any positive integer $ g < n$, there exists at least
one Hamming ball of radius $ n-g $, which contains at least ${n
\choose g}/q^{g-k} $ many codewords.  This bound matches the intuition
well, consider an imaginary algorithm as follows: randomly select $g$
points from the $n$ input points, and use polynomial interpolation to
get a polynomial of degree at most $g-1$ which passes these $g$
points. Then with probability $1/q^{g-k}$, the result polynomial has
degree $k -1 $. The sample space has size $n \choose g$. Thus
heuristically, the number of codewords in Hamming balls of radius
$n-g$ is at least $ {n \choose g}/q^{g-k} $ on the average.  In the
same paper, Justesen and Hoholdt also gave an upper bound for the
radius of the Hamming balls containing a constant or less number of
codewords. 

If we gradually increase $g$, starting from $k$, then
${n \choose g}/q^{g-k} $ will fall below 1 at some point.
However,   $g$ is still very far away from $\sqrt{nk}$.
Let $\hat{g}(n,k,q)$ be the smallest integer such that ${n \choose g}/q^{g-k} $
is less than $1$.
The following lemma shows that there is a gap between
  $\hat{g}(n,k,q)$ and $  \sqrt{nk} $.

\begin{lemma}\label{limitofGS}
\begin{enumerate}
\item For positive integers
$k < g < n $, if $g> \sqrt{nk}$, then
$ n^{ g-k} > {n \choose g}$
(which implies that $ q^{ g-k} > {n \choose g}$).
\item  For any constant $0 < c_1 < 1/2$
and fixed $k/n$, if $g = k + c_1 (n-k)$, then
${n \choose g}/n^{g-k} \leq 2^{-c_2 n} $ for some positive constant $c_2$.
\end{enumerate}
\end{lemma}

In fact, for a fixed rate ($k/n$) and $q = \Theta(n)$,
$\hat{g}(n,k,q) = k + \Theta ({ n \over \log n})$.  We prove that
if the
list decoding of the $[n,k]_q$ Reed-Solomon code is feasible when
radius is $n-\hat{g}(n,k,q)$, then the discrete logarithm over
$\F_{q^{\hat{g}(n,k,q)-k}}$ is easy.  In the other words,
we prove that the list decoding is not feasible for 
radius $n-\hat{g}(n,k,q)$ or larger, assuming that
the discrete logarithm over $\F_{q^{\hat{g}(n,k,q)-k}}$ is hard.
Note that it does not rule out
the possibility
that there are only polynomially many codewords in all Hamming balls
of radius  $ n -
\hat{g}(n,k,q) $, even assuming that intractability of
the discrete logarithm over $\F_{q^{\hat{g}(n,k,q ) -k}}$.

\begin{theorem}\label{MainI}
If there exists an algorithm solving the
list decoding problem of radius $n - \hat{g}(n,k,q) $
for the Reed-Solomon code $[n,k]_q$ in time $q^{O(1)}$,
then discrete logarithm over finite field $\F_{q^{\hat{g}(n,k,q) - k}}$
can be  computed in time $q^{O(1)} $.
\end{theorem}

When the list decoding problem is hard for certain radius, or a
Hamming ball contains too many codewords for us to enumerate all of
them, we can turn our attention to designing an efficient {\em bounded
distance decoding} algorithm, which only need to output one of
codewords in the ball, or output the empty set in case that the ball
does not contain any codeword.  However, we prove that the bounded distance
decoding is hard as well.

\begin{theorem}\label{MainII}

Let $q$ be a prime power and $h$ be a positive integer
satisfying $ q >  (h+2)^4 $. 
If the bounded distance decoding problem of radius $q- 4h - 4$
for the Reed-Solomon code $[q, 3h + 4]_q$ can be solved in
time $q^{O(1)} $,
the discrete logarithm problem over
$\F_{q^{h}}$ can be solved in time
$q^{O(1)} $.
\end{theorem}


To prove the theorem, we naturally come across the following question:
In a finite field $\F_{q^h}$, for any $\alpha$ such that $\F_{q^h} =
\F_q [\alpha]$, can $\F_q + \alpha$ generate the multiplicative group
$(\F_{q^h})^*$?  This interesting problem has a lot of applications in
graph theory, and it has been studied by several number theorists.
Chung \cite{chung89} proved that if $q > (h-1)^2$, then $(\F_{q^h})^*$
is generated by $\F_q + \alpha$.  Wan \cite{Wan97} showed a negative
result that if $q^h - 1$ has a divisor $d>1$ and $ h \geq 2 (q \log_q
d + \log_q (q+1))$, then $(\F_{q^h})^*$ is not generated by $\F_q +
\alpha$ for some $\alpha$.  Katz \cite{Katz90} applied the Lang-Weil
method, and showed that for every $h\geq 2$ there exists a constant
$B(h)$ such that for any finite field $\F_q$ with $q\geq B(h)$, any
element in $(\F_{q^h})^*$ can be written as a product of exactly $ n =
h +2$ distinct elements from $\F_q + \alpha$.  Clearly $B(h)$ has to
be an exponential function.  In this paper, we obtain a generalization
of the Katz's theorem, in which we use a bigger $n$ and manage to
decrease $B(h)$ to a polynomial function. For details, see
Section~\ref{BDD}

It is generally believed that the list decoding problem
and the bounded distance decoding for Reed-Solomon codes
are computationally
hard if the number of errors is greater than  $ n - \sqrt{nk}$ and
less than $n-k$.
This problem is even used as a hard problem to build
public key cryptosystems and pseudorandom generators \cite{KiayiasYu02}.
A similar problem, noisy polynomial
interpolation \cite{BleichenbacherNg00},
was proved to be vulnerable to the attack of
lattice reduction techniques, hence  is easier than originally thought.
This raises concerns  on the hardness of polynomial
reconstruction problem.
Our results confirm the belief that polynomial reconstruction
problem is hard, under a well-studied  hardness assumption
in number theory,
hence provide a firm foundation for
many protocols based on the problem.

This paper is organized as follows.
In Section~\ref{LemmalimitofGS}, we prove Lemma~\ref{limitofGS}.
In Section~\ref{maintheorem}, we sketch the proof of Theorem~\ref{MainI} and
Theorem~\ref{MainII}. In Section~\ref{size}, we show an
interesting duality between the  size of a group generated by linear factors,
and the list size in Hamming balls of Reed-Solomon codes.

\section{Proof of Lemma~\ref{limitofGS}}\label{LemmalimitofGS}

In this section, we prove Lemma~\ref{limitofGS} by
showing the following statement.

\begin{theorem}
There is no positive integral solution  for
\begin{eqnarray}
{{n \choose g} } & > & n^h\\
g &> & \sqrt{n (g - h)}.
\end{eqnarray}
\end{theorem}

We first obtain a finite range for $h, g$ and $n$.

\begin{lemma}
If $(n,g,h)$ is a positive integral solution, then $ h < 88. $
\end{lemma}

\begin{proof}
Denote $g/h$ by $\alpha$ and $n/h$ by $\beta$.
From $g >  \sqrt{n (g - h)}$, we have
$\alpha > \sqrt{\beta ( \alpha - 1)} $.
Hence $\alpha < \beta < \alpha + 1 + {1 \over \alpha - 1}. $

Recall that for any positive integer $i$, $ \sqrt{2\pi i}(i/e)^i \leq
i! \leq \sqrt{2\pi i}(i/e)^i ( 1 + {1 \over 12i-1}).  $

${n \choose g} = { \beta h \choose \alpha h}
 \leq ({ \beta^\beta \over \alpha^\alpha (\beta - \alpha)^{\beta - \alpha
 }})^h $.

Thus ${\beta^\beta \over \alpha^\alpha (\beta - \alpha)^{\beta - \alpha
 }} \geq \beta h $, which implies

$$ h \leq {\beta^{\beta-1} \over \alpha^\alpha (\beta - \alpha)^{\beta - \alpha
 }}.$$

Recall some facts:

\begin{enumerate}
\item For $x >0$, $x^x$ takes the minimum value $0.6922..$ at $x = e^{-1} =
0.36787944... $.
\item For $x>0$, $ 1 \leq (1 + { 1 \over x})^x \leq e = 2.7182818284...$
\end{enumerate}

If $\alpha \geq 2$, then $\beta - \alpha \leq 1 + { 1 \over \alpha -
  1} \leq 2$. We have
\begin{eqnarray*}
h  &\leq& {{1.45 \beta^{\beta - 1} \over \alpha^\alpha }}\\
   &\leq& {{1.45 ( 1 + \alpha +
 {1 \over \alpha -1})}^{( \alpha +
  {1 \over \alpha -1})}  \over \alpha^\alpha}\\
   &\leq& 1.45 {( 1+ \alpha + {1 \over \alpha - 1})^{( {1
  \over \alpha - 1})} ( 1 + { 1 \over \alpha}
   + { 1 \over \alpha ( \alpha-1)})^\alpha}\\
   &\leq& 1.45 * 4 * e * 2  < 32.
\end{eqnarray*}

If $\alpha < 2$, $ h \leq {1.45 \beta^{\beta-1} \over
(\beta-\alpha)^{\beta-\alpha}}$.
There are two cases. If $\beta \leq 3$, then
$$h \leq 1.45^2 * 9 < 19.  $$
If $\beta > 3$, then
\begin{eqnarray*}
h &\leq& 1.45 ({\beta \over \beta -
  \alpha})^{\beta - 1} (\beta-\alpha)^{\alpha - 1}\\
 &\leq& 1.45 ({\beta \over
  \beta - 2})^{\beta-1} (1 + {1 \over \alpha -1})^{\alpha - 1}\\
 &\leq& 1.45 * e^3 * 3 < 88.
\end{eqnarray*}

\end{proof}

\begin{corollary}
$\alpha \geq 88/87$ and $ \beta - \alpha < 88$.
\end{corollary}

Note that if $\alpha < 89$, then $\beta < 178$.
If $\alpha \geq 89$, then $\beta -\alpha \leq 1 + 1/88$,
but $n-g = (\beta-\alpha)h $ is an integer,
and $h \leq 87$, so $\beta - \alpha \leq 1$.
So if $n > 2h$, (1) can not hold.

\begin{proof}
Now we can finish proving the main
theorem of this section, by  exhaustively
searching for the solutions in the finite range that
$ h < 88, n < 178 * 88 = 15664$ and $h < g < n$ in a computer.
\end{proof}

Similarly we can show that for any constant $c$,
the inequalities
\begin{eqnarray}
{{n \choose g} } &\geq& n^{h-c}\\
g & > & \sqrt{t (g - h)}
\end{eqnarray}
have only finite many positive integral solutions.

Denote $n \over g -k$ by $\gamma $ and $g \over g-k$ by $\delta $.
To prove the second part of the lemma, it suffices to see that
$ {n\choose g} = { \gamma (g-k) \choose { \delta (g-k) }}  \leq c_2^{g-k}  $
for some constant $c_2$ only depending on $\alpha$ and $\beta$.

\section{The Decoding Problem of Reed-Solomon Codes
and the Discrete Logarithm over Finite Fields}\label{maintheorem}

Let $q$ be a prime power and let $\F_q$ be the finite field with
$q$ elements. Let $S$ be a subset of $\F_q$ of $n$ elements.
For a positive integer $g \leq n$, consider
$$ S_g = \{ A | A\subseteq S, |A| = g   \}. $$
For any $A \in S_g$, denote $\prod_{a \in A} (x-a)$ by $P_A (x)$.
Let $h(x)$ be an irreducible monic polynomial over $\F_q$ of degree $h
< g $.
Define a map $\psi: S_g \rightarrow \F_q[x]/(h(x))$ by
$$ \psi(A) = P_A(x) \pmod{h(x)}.    $$
For any $f(x)$ in $\F_q[x]/(h(x))$, if $\psi^{-1} ( f(x))$
is not empty, then there exists at least one polynomial $t(x)$
and one $A\in S_n$
such that
$f(x) + t(x) h(x) = P_A (x)$.
For any $a\in A$, $P_A (a) = 0$, $t(a) = - f(a)/h(a)$. Hence
there are at least $g$ elements  in $S $ which are
the roots of $f(x) + t(x) h(x) =0$, and the curve $y = t(x)$  passes
at least $g$ points in the following set of
$n$ points:
$$ \{ (a, - f(a)/h(a))| a\in S \}.  $$
According to Pigeonhole principle, there must exist a polynomial
$\hat{f} (x)$ such that $ | \psi^{-1} ( \hat{f}(x)) | \geq
|S_g |/|\F_q[x]/(h(x)) | = { {n \choose g } \over q^{h} }. $
Note that  $t(x)$ has degree
$g-h$ and  leading coefficient $1$.
For any polynomial $f \in \F_q [x]$  of degree at most $h-1$,
let $T_{f(x)}$ be the set of polynomial $t(x)$ of degree $g-h$
such that $f(x) + t(x)h(x) = P_A(x)$ for some $A\in S_g$,
and let $C_{f(x)}$ be the set of codewords within distance of $n-g$ to 
the received
word $( - f(a)/h(a) - a^{g-h} )_{a\in S} $ in Reed-Solomon code 
$[n,g-h]_q$.
{\em There is a one-to-one
correspondence between $T_{f(x)}$ and  $C_{f(x)}$,}
by sending any $t(x) \in T_{f(x)}$ to $ (t(a) - a^{g-h})_{a\in S} $.

Suppose that we  know $f(x)$ and $h(x)$, but not $A$,
are we still able to find $t(x)$? This is just a list
decoding problem of Reed-Soloman code $[n, g-h]_q$. Once we have
a list of $t(x)$, we can find $A$ by factoring $f(x) + t(x) h(x)$.
This provides  a general framework for the following proofs.

\subsection{The proof of Theorem~\ref{MainI} }

Given a Reed-Solomon code $[n,k]_q$, let $h = \hat{g}(n,k,q)-k $.
Recall that $\hat{g}(n,k,q)$ is the
smallest integer such that ${n \choose g}/q^{g-k} $
is less than $1$, and $h$ is the degree of an irreducible polynomial $h(x)$.
We show that there is an efficient algorithm to solve
the discrete logarithm over $\F_{q^{\hat{g}(n,k,q)-k}} = \F_q[x]/(h(x))$
if there is efficient list decoding algorithm for the Reed-Solomon
code $[n,k]_q$ with radius $n-\hat{g}(n,k,q)$.
Let $\alpha = x \pmod{h(x)}$.
Suppose that we are given the base $b(\alpha)$ and we need to
find out the discrete logarithm of $t(\alpha)$ with respect to the base,
where $b$ and $t$ are polynomials over $\F_q$ of degree at most
$h - 1$.
That there is an efficient list decoding algorithm implies:
\begin{enumerate}
\item There are only polynomially many codewords in any Hamming ball
of radius $n-\hat{g}(n,k,q)$, which in turn implies that $ | \psi^{-1}
(f) | \leq q^c $ for any $f\in \F_{q^h}$ and a constant $c$. Hence
$$ |\psi (S_{\hat{g}(n,k,q)}  )| \geq
{{ n \choose \hat{g}(n,k,q)   } \over q^c} =
\Theta (q^{\hat{g}(n,k,q) - k}/q^c )=\Theta (q^h/q^c).  $$
\item And they can be found in polynomial time.
\end{enumerate}
We use the index calculus algorithm with {\em factor bases}
$(\alpha + a)_{a \in S}$.
If we randomly select an integer $i$ between $0$ and
$q^{\hat{g}(n,k,q)-k} - 1$,
then with  probability bigger than $1/q^c$, $\psi^{-1}(b(\alpha)^i)$ is
not empty. Apply the list decoding algorithm,
we get  {\em relations}
$$ b(\alpha)^i = f(\alpha) = \prod_{a\in A_1 } (\alpha + a)
= \cdots
= \prod_{a\in A_l } (\alpha + a) $$ for some $A_1, A_2, \cdots, A_l \in
S_{\hat{g}(n,k,q)} $ where $l$ is the list size.
From the relations, we get  linear equations.

$$ i = \sum_{a \in A_1} \log_b (\alpha + a) = \cdots
= \sum_{a \in A_l } \log_b (\alpha + a)
 \pmod{q^{\hat{g}(n,k,q)-k} - 1  }$$

We repeat the above procedure.
Since $i$ is picked randomly, and $S_g$ is the sample space,
the probability that the new equation is linear independent to the
previous ones is very high at the beginning of the algorithm.
It would not take long time
before we get $n$ independent equations.
Solving the system of equations gives us
$ \log_b (\alpha + a)$ for all $ a \in \F_q  $.

In the last step, for a random $i$, we compute $b(\alpha)^i t(\alpha)$.
If $\psi^{-1} (b(\alpha)^i t(\alpha))$ is not empty, we can solve
$\log_b t$ immediately.
This proves the main theorem.

\subsection{The proof of Theorem~\ref{MainII} }
\label{BDD}

\begin{theorem}\label{product}
Let $q$ be a prime power and let $h$ be a positive integer.
    If $q \geq (h+2)^4$,
    then  every element in $\F_{q^h}^*$
    can be written as a product of exactly $4h+4$ distinct
    factors from $\{ \alpha + a | a\in F_q \}$,
    for any $\alpha$ such that $\F_q(\alpha) = \F_{q^h}$.

\end{theorem}

\begin{proof} We thank Chaohua Jia for helpful discussion on the
proof of this theorem.  Fix an $\alpha$ such that $\F_q(\alpha) =
\F_{q^h}$. For $\beta\in \F_{q^h}^*$, let $N_k(\beta)$ denote the
number of solutions of the equation
$$\beta = \prod_{i=1}^k (\alpha +a_i), \ a_i \in \F_q,$$
where the $a_i$'s are distinct. We need to show that for $k=4h+4$,
the number $N_k(\beta)$ is always positive if $q \geq (h+2)^4$.

Let $G$ be the character group of the multiplicative group
$\F_{q^h}^*$, which is a cyclic group of order $q^h-1$. A simple
inclusion-exclusion argument shows that
$$N_k(\beta) \geq {1\over q^h-1} (\sum_{a_i\in \F_q, 1\leq i\leq
k} - \sum_{1\leq i_1<i_2\leq k}\sum_{a_i\in \F_q,
a_{i_1}=a_{i_2}}) \sum_{\chi\in G} \chi^{-1}(\beta)
\chi(\prod_{i=1}^k (\alpha+a_i)).$$ For non-trivial $\chi$, one
has the well-known Weil estimate
$$|\sum_{a\in \F_q}\chi(\alpha+a)| \leq (h-1)\sqrt{q}.$$
We deduce that
$$N_k(\beta) \geq {q^k -{k\choose 2}q^{k-1} \over q^h-1}
-(1+{k\choose 2})(h-1)^k q^{k/2}.$$ In order for $N_k(\beta)>0$,
it suffices to have the inequality
$$(q-{k\choose 2})q^{k/2 -1-h} > (1+{k\choose 2}(h-1)^k.$$
This inequality is clearly satisfied if both $q>2{k\choose
2}+1=k(k-1)+1$ and $q^{k/2 -1-h} >(h-1)^k$. These two inequalities
are satisfied if we take $k=4h+4$ and $q\geq (h+2)^4$. The theorem
is proved.

\end{proof}

Now we are ready to prove Theorem~\ref{MainII}

\begin{proof}
Let $ h(x)$ be an irreducible polynomial
over $\F_q$ of degree $h$.
Then $\F_{q^h} = \F_q [x]/(h(x)).$
Denote $x \pmod{h(x)}$ as $\alpha$.
Suppose we need to solve the discrete logarithm of $t(\alpha)$ base
$b (\alpha)$ in $\F_{q^h}$, where $b$ and $t$ are polynomials
of degree at most $h-1$. We let $ S = \F_q $.
$$ (\F_q)_{4h+4} = \{ A | A\subseteq \F_q, |A| = 4h+4   \}. $$

First we randomly select an integers $i$ between
$0$ and $q^h-1$. Compute $b(\alpha)^i$,
and let $f(\alpha) $ be the result where $f(x)$ is
a polynomial of degree at most $h-1$. Now run
the bounded distance decoding algorithm on
the  Reed-Solomon code $[q, 3h + 4]_q$ with the
point set $ \{(a, - f(a)/h(a) - a^{3h+4 }) |  a\in \F_q \} $ and
the distance bound $ q - 4h-4$.
Then according to
Theorem~\ref{product}, the answer is not the empty set.
Let the answer be $t(x) - x^{ 3h+4 }$. The polynomial
$t(x)$  has degree $ 3h+4$,
and agrees with
$\{ (x, - f(x)/h(x)) | x \in \F_q \}$ at  $ 4h+4 $ many points or more.
The polynomial $f(x) + t(x) h(x)$ has degree at most $ 4h+4$,
but has at least $4h+4$ many distinct zeros, thus
it will be completely splitted
as a product of
linear factors.
Let $f(x) + t(x) h(x) = \prod_{ a\in A } ( x  + a)$
for some $ A \in (\F_q)_{4h+4} $.
Write it in another way,
$$ b^i = \prod_{ a\in A } ( \alpha  + a).$$
We  get
$$ i = \sum_{a \in A} \log_g ( \alpha + a) \pmod{q^h-1}.$$

However, we may not be able to solve $\log_g (\alpha + a)$
for all $a\in \F_q$,
since the latter relations may be linearly dependent on
the former relations. This is the case, for instance,  when all the $A_i$'s
come from a subset of $\F_q$.
After we detect that,
we start to compute
$t(\alpha) b(\alpha)^x$, and find its
representation of product of linear factors.
Any linear dependence will give us the discrete
logarithm of $t(\alpha)$ base $ b(\alpha)$.
\end{proof}

%
%
%

\section{Group Size and List Size}\label{size}

Let $q$ be a prime power, and $S$ be a subset of $\F_q$ of $n$
elements, where $n $ is very small compared to $q$.  Let $\alpha $
be an element in $\F_{q^h}$ such that $\F_q[\alpha] = \F_{q^h}$. What
is the order of the subgroup generated by $\alpha + S$ for some
$S\subseteq \F_q$ ?  This question has an important application in
analyzing the performance of the AKS primality testing
algorithm \cite{AgrawalKa02}. 
Experimental data suggests that the order is greater than $ q^{h/c}$ for some
absolute constant $c$ for $|S| \geq h\log q$.
If we can prove it, the space  complexity of the AKS algorithm can be
cut by a factor of $\log p$ ($p$ is the input prime whose
primality certificate is sought),
which will make (the random variants of ) the algorithm
comparable to the primality proving algorithm used in practice.
However, the best known lower bound  is $(c|S|/h)^h$ for some
absolute constant $c$ \cite{Voloch02}.
We  discover an interesting
duality  between the group size and the list size in Hamming
balls of certain radius.

%

\begin{theorem}
Let $k,n$ be positive integers and $q$ be a prime power.
One of the following statements must be true.
\begin{enumerate}
\item For any constant $c_1$, there exists
a Reed-Solomon code $[n,k]_q$ ($n/3 < k < n/2$),
and a
Hamming ball of radius $n- \hat{g}(n,k,q)$ containing
more than $c_1 1.9^n $  codewords.
\item Let $s = \log q$, the
group generated by $\alpha + S $,
has cardinality at least $ q^{h/c_2} $ for some absolute constant $c_2$,
where  $S \subseteq \F_q $ and $ |S| = s\log q$.
\end{enumerate}
\end{theorem}

To prove the first statement would solve an
important open problem in the Reed-Solomon codes.
To prove the second statement would give us
a primality proving algorithm much more efficient
in term of  space complexity than the original AKS and
its random variants, hence make
the AKS algorithm not only theoretical
interesting, but also practical important.
However, at this stage we cannot figure out
which one is true. What we can prove, however, is that
one of them must be true. Note that it is also possible
that both of the statements are true.


\begin{proof}
Let $s = \log q$,  $k = sh/2 - h$ and $n = sh$.
So the rate $k/n$ is very close
to $1/2$ as $s$ gets large, and $\hat{g}(n,k,q) = sh/2$.
Assume the first statement is wrong,
this means that there exists a constant $c_3$ such that
for any Reed-Solomon code $[n,k]_q$ with $n/3 < k < n/2 $,
the number of codewords in any Hamming ball of radius $ n - \hat{g}(n,k,q) $
is less than $c_3 1.9^n$.
The number of  balls containing at least one codeword 
with that radius and center point
at $ ( - f(a)/h(a) - a^{k})_{a \in S \in\F_q}  $,
where $f\in \F_q[x]$ has degree less than
$ h $
is greater than
$$  q^h / (c_3 1.9^n) = q^{ h - n\log 1.9 /\log q}/c_3 \geq q^{h/c}, $$
which is a low bounded of the size of the group generated
by $\alpha + S$.
\end{proof}

\section{Concluding Remarks}

Interesting open questions include whether the decoding problem of
Reed-Solomon code is equivalent to or harder than the
discrete logarithm over finite fields, and
whether there exists  a polynomial time quantum algorithm to solve the
decoding problem of Reed-Solomon code.


\bibliographystyle{plain}
\bibliography{crypto}

\end{document}